\documentclass[a4paper,12pt]{article}
\usepackage{amsmath,amsfonts,amssymb,amsthm}  
\newtheorem{theo}{Theorem}[section]

\newcommand{\la}{\lambda}
\newcommand{\ka}{\theta}

\begin{document}
\title{An analytic formula for Macdonald 
polynomials \\
{\small Une formule analytique pour les polyn\^omes de Macdonald}}
\author{Michel {\sc Lassalle} and Michael {\sc Schlosser}}
\date{5 mai 2003}
\maketitle
\begin{abstract}
We give the explicit analytic development of any Jack or Macdonald
polynomial in terms of elementary
(resp.\ modified complete) symmetric functions. These two
developments are obtained by inverting the Pieri formula. \\

\begin{center}
\textbf{R\'esum\'e}   
\end{center}
Nous donnons le d\'eveloppement analytique explicite de tout
polyn\^ome de Jack ou de Macdonald sur les fonctions
sym\'etriques \'el\'ementaires (resp.\ compl\`etes modifi\'ees).
Nous obtenons ces deux d\'eveloppements par inversion de la formule
de Pieri.
\end{abstract}
                           
\begin{center}
\textbf{Version fran\c{c}aise abr\'eg\'ee}    
\end{center}

Au milieu des ann\'ees cinquante, Hua introduisait les
polyn\^omes zonaux et posait le probl\`eme d'en obtenir un
d\'eveloppement analytique explicite~\cite{H}. En d\'epit de nombreuses
recherches, cette question est demeur\'ee ouverte.

Elle est d\'esormais formul\'ee dans le cadre plus g\'en\'eral des
polyn\^omes de Macdonald. Les polyn\^omes zonaux sont en effet
un cas particulier des polyn\^omes de Jack, qui sont eux-m\^emes
un cas limite des polyn\^omes de Macdonald.

Ces polyn\^omes sont index\'es par
les partitions d'entiers. Ils forment
une base de l'alg\`ebre des fonctions sym\'etriques \`a
coefficients rationnels en deux param\`etres $q,t$, et
g\'en\'eralisent notamment les fonctions de Schur, les
polyn\^omes de Hall--Littlewood, et les polyn\^omes de Jack.

Les polyn\^omes de Macdonald ont \'et\'e d\'etermin\'es par Lapointe,
Lascoux et Morse \cite{LLM1}, qui les ont exprim\'es
comme d\'eterminants. Cependant les entr\'ees de ces 
d\'eterminants sont des quantit\'es combinatoires qui 
ne peuvent \^etre en g\'en\'eral exprim\'ees analytiquement.

On ne disposait donc jusqu'ici 
d'aucune formule analytique explicite 
pour les polyn\^omes de Jack et de Macdonald. On ignorait notamment leur 
d\'eveloppement sur les bases classiques, 
sauf dans deux cas particuliers : lorsque la partition indexante est de 
longueur 
deux \cite{JJ}, ou lorsqu'elle est de longueur trois \cite{La}, ainsi que 
dans 
les situations duales o\`u les parts sont au plus 
\'egales \`a $2$ ou $3$.

Le but de cette Note est de pr\'esenter une solution g\'en\'erale
\`a ce probl\`eme. Pour tout polyn\^ome de Macdonald nous
obtenons deux d\'eveloppements
analytiques explicites. Le premier se fait sur les
fonctions sym\'etriques \'el\'ementaires. Le
second sur les fonctions sym\'etriques ``compl\`etes
modifi\'ees'', dont le d\'eveloppement sur toute base classique
est connu \cite{La1}.

\medskip

\section{Introduction}

In the fifties, Hua posed the problem of finding
an explicit analytic formula for the ``zonal polynomials'' he
had just introduced~\cite{H}. In spite of many efforts this
question has remained open.

Hua's problem is now better understood in the 
more general framework of Macdonald
polynomials. Zonal polynomials are indeed a special case of
Jack polynomials, which in turn are obtained from
Macdonald polynomials by taking a particular limit.

Macdonald polynomials are indexed by
partitions. They form a basis of the algebra
of symmetric functions with rational coefficients in two
parameters $q,t$, and generalize Schur functions,
Hall--Littlewood polynomials and Jack polynomials.

Lapointe, Lascoux and Morse~\cite{LLM1} gave a determinantal expression
for Macdonald polynomials.
However their method does not lead to an explicit analytic
formula, since the entries of these determinants are
combinatorial quantities which in general cannot be analytically
explicited.

Thus Hua's problem kept open for Macdonald polynomials.
Their analytic expansion was only known when the length of the
indexing partition is two~\cite{JJ}, or three~\cite{La},
and in the dual cases corresponding to partitions with
parts at most equal to $2$ or $3$.

The aim of this Note is to present a general solution
to this problem. We give two explicit analytic developments
for any Macdonald polynomial.
One of them is made in terms of elementary
symmetric functions. The other one is made in terms of
``modified complete'' symmetric functions, which have 
themselves a known development in terms of 
any classical basis~\cite{La1}.

Our method relies on the inversion of the ``Pieri matrix''.
This is an infinite transition matrix which has been
analytically explicited by Macdonald~\cite{Ma}. We compute its inverse
by the operator method of Krattenthaler~\cite{K1,K2}, such as
adapted to the multivariate case by the second
author~\cite{S1}.

Detailed proofs will be given in a forthcoming paper.

\section{Macdonald Polynomials}

The standard reference for Macdonald polynomials is Chapter 6
of \cite{Ma}.

Let $X=\{x_1,x_2,x_3,\ldots\}$ be an infinite set
of indeterminates (an alphabet) and $\mathcal{S}$ the
corresponding algebra of symmetric functions
with coefficients in $\mathsf{Q}$.
Elementary symmetric functions $e_{k}(X)$,
complete symmetric functions $h_{k}(X)$ and
power sums $p_{k}(X)$ form three algebraic bases
of $\mathcal{S}$.

Let $q,t$ be two indeterminates. We define ${(a;q)}_{0}:=1$,
${(a;q)}_{k}:=\prod_{i=0}^{k-1}(1-aq^i)$, for $k\ge 1$, and
${(a;q)}_{\infty}:=\prod_{i\geq0}(1-aq^i).$

A partition $\la= (\la_1,...,\la_n)$
is a finite weakly decreasing
sequence of positive integers, called parts. The number
$n=l(\la)$ of parts is called the length of
$\la$. For any integer $i\geq1$,
$m_i(\la) := \textrm{card} \{j: \la_j  = i\}$
is the multiplicity of $i$ in $\la$.  Clearly
$l(\la)=\sum_{i\ge1} m_i(\la)$. We shall also write
 $\la= (1^{m_1},2^{m_2},3^{m_3},\ldots)$.

Let $\mathsf{Q}[q,t]$ be the algebra of rational functions  in
$q,t$, and $\mathsf{Sym}=\mathcal{S}\otimes\mathsf{Q}[q,t]$
the algebra of symmetric functions
with coefficients in $\mathsf{Q}[q,t]$. This algebra is
endowed with some scalar product for which Macdonald polynomials
$P_{\la}(X;q,t)$ form an orthogonal basis.
Let $Q_{\la}(X;q,t)= b_{\la}(q,t) P_{\la}(X;q,t)$ be the dual
basis, with $b_{\la}(q,t)$ specified in \cite{Ma}, p.~339, Equ.~(6.19).

For any $k\ge 0$ the ``modified complete'' symmetric function
$g_{k}(X;q,t)$ is defined by the generating series
\begin{equation*}
\prod_{i \ge 1}
\frac{{(tux_i;q)}_{\infty}}{{(ux_i;q)}_{\infty}}
=\sum_{k\ge0} u^k g_{k}(X;q,t).
\end{equation*}
Then it is known (\cite{Ma}, p.~329, Equ.~(5.5)) that
the Macdonald polynomial associated with a row partition
$(k)$ is given by  $Q_{(k)}(q,t)= g_{k}(q,t)$.

The symmetric functions $g_{k}(q,t)$ form an algebraic basis of
$\mathsf{Sym}$. They may be expanded in terms of any classical
basis. This
development is explicitly given in \cite{Ma} (p.~311 and 314)
in terms of power sums and monomial symmetric functions,
and in \cite{La1} (Sec.~10, p.~237) for other classical bases.

The parameters $q,t$ being kept fixed, we shall write $P_{\mu}$
for $P_{\mu}(q,t)$ (resp.\ $Q_{\mu}$ for $Q_{\mu}(q,t)$).

\section{Pieri formula}

Let $u_1,\ldots,u_n$ be $n$ indeterminates and $\mathsf{N}$ the set of
nonnegative integers. For $\ka =(\ka_1,\ldots,\ka_n)\in \mathsf{N}^{n}$
let $|\ka|=\sum_{i=1}^{n} \ka_i$ and define
\begin{equation*}
d_{\ka_1,\ldots,\ka_n} (u_1,\ldots,u_n):=
\prod_{k=1}^n \frac{{(t;q)}_{\ka_k}}{{(q;q)}_{\ka_k}}\,
\frac{{(q^{|\ka|+1}u_k;q)}_{\ka_k}}{{(q^{|\ka|}tu_k;q)}_{\ka_k}}
\prod_{1\le i < j \le n}
\frac{{(tu_i/u_j;q)}_{\ka_i}}{{(qu_i/u_j;q)}_{\ka_i}}\,
\frac{{(q^{-\ka_j+1}u_i/tu_j;q)}_{\ka_i}}{{(q^{-\ka_j}u_i/u_j;q)}_{\ka_i}}.
\end{equation*}

Macdonald polynomials satisfy a Pieri
formula (\cite{Ma}, p.~340, Equ.~(6.24) (ii)), which may be
analytically explicited as follows (\cite{Ma}, p.~342, Example 2(b)).

\begin{theo}[Macdonald]\label{theopieri}
Let  $\la=(\la_1,...,\la_n)$ be an arbitrary partition with length
$n$ and $\la_{n+1} \in \mathsf{N}$. For all $1 \le k \le n$ define
$u_k:=q^{\la_k-\la_{n+1}}t^{n-k}$. We have
\begin{equation*}
Q_{(\la_1,\ldots,\la_n)} \: Q_{(\la_{n+1})}= \sum_{\ka\in \mathsf{N}^n}
d_{\ka_1,\ldots,\ka_n} (u_1,\ldots,u_n)\:
Q_{(\la_1+\ka_1,\ldots,\la_n+\ka_n,\la_{n+1}-|\ka|)}.
\end{equation*}
\end{theo}

This Pieri formula defines an infinite transition matrix.
Indeed let $\mathsf{Sym}(n+1)$ denote the algebra of symmetric
polynomials
in $n+1$ independent variables with coefficients in $\mathsf{Q}[q,t]$.
Then the Macdonald polynomials
$\{Q_\la, \,l(\la) \le n+1\}$ form a basis of
$\mathsf{Sym}(n+1)$, and so do the
products $\{Q_{\mu}Q_{(k)}, \,l(\mu) \le n, \,k \ge 0\}$.

\section{Main results}\label{sectmain}

Let $u=(u_1,\ldots,u_n)$ and $v=(v_1,\ldots,v_n)$ be $2n$ indeterminates.
We denote by $\Delta(v)$ the Vandermonde determinant
$\prod_{1\le i < j \le n} (v_i-v_j)$.
For $a,b$ two indeterminates and
$\ka =(\ka_1,\ldots,\ka_n)\in \mathsf{N}^{n}$ we define
\begin{multline*}
c^{(a,b)}_{\ka} (u;v):= \prod_{k=1}^n
b^{\ka_k} \,\frac{(a/b;a)_{\ka_k}}{(a;a)_{\ka_k}}\,
\frac{(au_k;a)_{\ka_k}}{(abu_k;a)_{\ka_k}}\,
\prod_{1\le i < j \le n}
\frac{{(au_i/bu_j;a)}_{\ka_i}}{{(au_i/u_j;a)}_{\ka_i}} \,
\frac{{(a^{-\ka_j}bu_i/u_j;a)}_{\ka_i}}{{(a^{-\ka_j}u_i/u_j;a)}_{\ka_i}}\\
\times\frac{1}{\Delta(v)} \,
\operatorname{det}\!
{\left[v_i^{n-j}
\left(1-b^{j-1} \frac{1-bv_i}{1-v_i}
\prod_{k=1}^n \frac{u_k-v_i}{bu_k-v_i}\right)\right]}_{1\le i,j \le n}.
\end{multline*}

We have the following result.

\begin{theo}\label{theomain}
Let  $\la=(\la_1,...,\la_{n+1})$ be an arbitrary partition with length
$n+1$. For all $1\le k \le n$ define
$u_k:=q^{\la_k-\la_{n+1}}t^{n-k}$. We have
\begin{equation}\label{theomainid}
Q_{(\la_1,\ldots,\la_{n+1})}= \sum_{\ka\in\mathsf{N}^n}
c^{(q,t)}_{\ka_1,\ldots,\ka_n} (u_1,\ldots,u_n;
u_1q^{\ka_1},\ldots,u_nq^{\ka_n})\:
Q_{(\la_{n+1}-|\ka|)} \:
Q_{(\la_1+\ka_1,\ldots,\la_n+\ka_n)}.
\end{equation}
\end{theo}

The reader will easily check that for $n=1$ we obtain the
result previously given in \cite{JJ}. For
$n=2$ we recover the formula announced in a previous note
by the first author \cite{La}.

\begin{proof}[Sketch of proof of Theorem~\ref{theomain}]
We invert the infinite transition matrix defined by the Pieri formula
in Theorem~\ref{theopieri}. We accomplish this by applying
Krattenthaler's~\cite{K1,K2} operator method for
proving matrix inversions as adapted to the multivariate case
by the second author~\cite{S1}.

Let $\mathsf Z$ denote the set of integers and let
$\beta=(\beta_1,\dots,\beta_n)$,
$\kappa=(\kappa_1,\dots,\kappa_n)$,
$\gamma=(\gamma_1,\dots,\gamma_n)\in\mathsf Z^n$.
Defining
\begin{equation}\label{deff}
f_{\beta\kappa}:=c_{\beta_1-\kappa_1,\ldots,\beta_n-\kappa_n}^{(q,t)}
\big(u_1q^{\kappa_1+|\kappa|},\ldots,u_nq^{\kappa_n+|\kappa|};
u_1q^{\beta_1+|\kappa|},\ldots,u_nq^{\beta_n+|\kappa|}\big),
\end{equation}
and
\begin{equation}\label{defg}
g_{\kappa\gamma}:=d_{\kappa_1-\gamma_1,\ldots,\kappa_n-\gamma_n}
\big(u_1q^{\gamma_1+|\gamma|},\ldots,u_nq^{\gamma_n+|\gamma|}\big),
\end{equation}
our method yields that the infinite lower-triangular multidimensional 
matrices
$(f_{\beta\kappa})_{\beta,\kappa\in\mathsf Z^n}$ and
$(g_{\kappa\gamma})_{\kappa,\gamma\in\mathsf Z^n}$
are {\em inverses} of each other, i.~e., the
orthogonality relation
\begin{equation}\label{orthrel}
\sum_{\kappa\in\mathsf Z^n}f_{\beta\kappa}g_{\kappa\gamma}=
\delta_{\beta\gamma},
\end{equation}
for all $\beta,\gamma\in\mathsf Z^n$, holds.

It is immediately clear from \eqref{orthrel} that if
$(w_{\kappa})_{\kappa\in\mathsf Z^n}$ and $(y_{\kappa})_{\kappa\in\mathsf 
Z^n}$
are two multivariate sequences of indeterminates,
then
\begin{equation}\label{invrelf}
\sum_{\beta\in\mathsf Z^n}f_{\beta\kappa}w_{\beta}=y_{\kappa}\qquad\qquad
\text{for all $\kappa$} \in\mathsf Z^n,
\end{equation}
if and only if
\begin{equation}\label{invrelg}
\sum_{\kappa\in\mathsf 
Z^n}g_{\kappa\gamma}y_{\kappa}=w_{\gamma}\qquad\qquad
\text{for all $\gamma$} \in\mathsf Z^n,
\end{equation}
subject to convergence.

Now if in Theorem~\ref{theopieri}, we replace $u_i$ by
$u_iq^{\gamma_i+|\gamma|}$, $\lambda_i$ by
$\lambda_i+\gamma_i$, for $i=1,\dots,n$, and $\lambda_{n+1}$ by
$\lambda_{n+1}-|\gamma|$,
then we see (after shifting the summation indices)
that \eqref{invrelg} holds for
$y_{\kappa}=Q_{(\la_1+\kappa_1,\dots,\la_n+\kappa_n,\la_{n+1}-|\kappa|)}$
and $w_{\gamma}=Q_{(\la_1+\gamma_1,\ldots,\la_n+\gamma_n)}
Q_{(\la_{n+1}-|\gamma|)}$, with $g_{\kappa\gamma}$ given as in 
\eqref{defg}.
Thus, with $f_{\beta\kappa}$ given as in \eqref{deff},
we immediately establish \eqref{invrelf} for above values
of $y_{\kappa}$ and $w_{\gamma}$.
Setting now $\kappa_i=0$ for $i=1,\dots,n$ gives \eqref{theomainid}.
\end{proof}

It is known (\cite{Ma}, p.~327) that there exists an automorphism
$\omega_{q,t}$ of $\mathsf{Sym}$ such that
\[\omega_{q,t}(Q_{\la}(q,t))=P_{\la^{'}}(t,q),\quad
\omega_{q,t}(g_{n}(q,t))=e_{n},\]
with $\la^{'}$ the partition conjugate to $\la$, whose parts
are given by $m_k(\la^{'})=\la_k-\la_{k+1}$.
Applying this automorphism to the previous theorem we obtain
the following equivalent result.
\begin{theo}\label{theodual}
Let  $\la=(1^{m_1},2^{m_2},\ldots,(n+1)^{m_{n+1}})$ be an arbitrary
partition consisting of parts at most equal to $n+1$.
For all $1\le k \le n$ define
$u_k:=q^{n-k}t^{\sum_{j=k}^n m_j}$. We have
\begin{multline*}
P_{(1^{m_1},2^{m_2},\ldots,(n+1)^{m_{n+1}})}=
\sum_{\ka\in\mathsf{N}^n}
c^{(t,q)}_{\ka_1,\ldots,\ka_n} (u_1,\ldots,u_n;
u_1t^{\ka_1},\ldots,u_nt^{\ka_n}) \\
\times e_{m_{n+1}-|\ka|} \,
P_{(1^{m_1+\ka_1-\ka_2},2^{m_2+\ka_2-\ka_3},\ldots,
{(n-1)}^{m_{n-1}+\ka_{n-1}-\ka_n},n^{m_n+m_{n+1}+\ka_n})}.
\end{multline*}
\end{theo}

\section{Analytic developments}\label{sectanal}

Theorem \ref{theomain} (resp.\ Theorem \ref{theodual}) immediately 
generates
the analytic development of any Macdonald polynomial in 
terms of the
symmetric functions $g_k$ (resp.\ the elementary symmetric functions
$e_k$), which form an algebraic basis of $\mathsf{Sym}$.

Indeed for any partition $\mu=(\mu_1,...,\mu_l)$ with length $l > 1$ let
us write
\[c_{\ka_1,\ldots,\ka_{l-1}}(\mu):=
c^{(q,t)}_{\ka_1,\ldots,\ka_{l-1}}
(u_1,\ldots,u_{l-1};
u_1q^{\ka_1},\ldots,u_{l-1}q^{\ka_{l-1}})\]
with $u_k:=q^{\mu_k-\mu_l}t^{l-k-1}$.
Let $\mathsf{M}^{(n)}$ be the set of lower triangular $n \times n$
matrices with nonnegative integers.
By an obvious iteration we deduce
the following analytic expansion of any Macdonald polynomial.
\begin{theo}
Let  $\la=(\la_1,...,\la_{n+1})$ be an arbitrary partition with length
$n+1$. For any $\ka=(\ka(i,j))_{i,j=1}^n \in 
\mathsf{M}^{(n)}$, let us consider
a sequence of $n$ partitions $\{\mu(\ka,k), 1\le k\le n\}$
where $\mu(\ka,k)$ has length $k+1$ and is defined by
\[\mu(\ka,k)_{i}=\la_{i}+\sum_{j=k+1}^{n} \ka(j,i)
\quad\quad (1 \le i \le k+1).\]
We have
\[Q_{\la}= \sum_{\ka\in \mathsf{M}^{(n)}}
\prod_{k=1}^n
c_{\ka(k,1) \ldots \ka(k,k)}(\mu(\ka,k)) \: \prod_{k=0}^n
g_{\la_{k+1}+\sum_{j=k+1}^{n} \ka(j,k+1)-\sum_{j=1}^{k} \ka(k,j)}.
\]
\end{theo}

For any partition $\mu=(1^{m_1},2^{m_2},\ldots,l^{m_l})$ 
having its parts
at most equal to $l$, let us write
\[C_{\ka_1,\ldots,\ka_{l-1}}(\mu):=
c^{(t,q)}_{\ka_1,\ldots,\ka_{l-1}}
(u_1,\ldots,u_{l-1};
u_1t^{\ka_1},\ldots,u_{l-1}t^{\ka_{l-1}})\]
with $u_k:=q^{l-k-1}t^{\sum_{j=k}^{l-1} m_j}$.
By duality we deduce the following analytic expansion of
Macdonald polynomials in terms of elementary symmetric functions.
\begin{theo}
Let  $\la=(1^{m_1},2^{m_2},\ldots,(n+1)^{m_{n+1}})$ be an arbitrary
partition consisting of parts at most equal to $n+1$.
For any $\ka \in \mathsf{M}^{(n)}$ let
us consider a sequence of $n$ partitions
$\{\mu(\ka,k), 1\le k\le n\}$
where $\mu(\ka,k)$ has parts at most equal to $k+1$ and is defined by
\begin{equation*}
\begin{split}
m_i(\mu(\ka,k))&=m_{i}+\sum_{j=k+1}^{n}
\left(\ka(j,i)- \ka(j,i+1)\right)
\quad\quad (1 \le i \le k)\\
m_{k+1}(\mu(\ka,k))&=\sum_{j=k+1}^{n+1} m_{j}+\sum_{j=k+1}^{n} \ka(j,k+1).
\end{split}
\end{equation*}
We have
\[P_{\la}= \sum_{\ka\in \mathsf{M}^{(n)}}
\prod_{k=1}^n
C_{\ka(k,1) \ldots \ka(k,k)}(\mu(\ka,k)) \: \prod_{k=0}^n
e_{\,\sum_{j=k+1}^{n+1}m_j+\sum_{j=k+1}^{n} \ka(j,k+1)-\sum_{j=1}^{k}
\ka(k,j)}.
\]
\end{theo}

\noindent
\textbf{Remark :}
The expression given by Theorem 5.1 can be written in 
terms of raising operators (\cite{Ma}, p.~9). 
Writing $g_{\mu} = \prod_{k \ge 1} g_{\mu_k}$, the 
raising operators $R_{ij}$ act as follows : 
$R_{ij}g_{\mu}=g_{\mu_{1}}\ldots g_{\mu_{i}+1}\ldots g_{\mu_{j}-1}\ldots$.
Then $\prod_{k=0}^n g_{\la_{k+1}+\sum_{j=k+1}^{n} 
\ka(j,k+1)-\sum_{j=1}^{k} \ka(k,j)}$ 
is exactly $\big(\prod_{1\le i<j \le n+1}R_{ij}^{\ka(j-1,i)}\big) \, 
g_{\la}$.

\section{Jack polynomials}

Jack polynomials are
the limit of Macdonald polynomials when $t\rightarrow 1$, with
$q=t^\alpha$
and $\alpha$ some positive real number (\cite{Ma}, p.~376). We note
$ P_{\la} = \lim_{t \rightarrow 1} P_{\la}(t^\alpha,t)$ and
$ Q_{\la} = \lim_{t \rightarrow 1} Q_{\la}(t^\alpha,t)$.

Let ${(u)}_k$ be the classical raising factorial defined by
$(u)_k=\prod_{i=1}^{k} (u+i-1)$.
Let $u=(u_1,\ldots,u_n)$ and  $v=(v_1,\ldots,v_n)$  be $2n$
indeterminates.
For $\ka =(\ka_1,\ldots,\ka_n)\in \mathsf{N}^{n}$ and any
indeterminate $a$ define
\begin{multline*}
c^{(a)}_{\ka} (u;v) := \prod_{k=1}^n
\frac{(1-a)_{\ka_k}}{\ka_k!} \,
\frac{(u_k+1)_{\ka_k}}{(u_k+1+a)_{\ka_k}}
\prod_{1\le i < j \le n}
\frac{{(u_i-u_j+1-a)}_{\ka_i}}{{(u_i-u_j+1)}_{\ka_i}} \,
\frac{{(u_i-u_j-\ka_j+a)}_{\ka_i}}{{(u_i-u_j-\ka_j)}_{\ka_i}} \\
\times \frac{1}{\Delta(v)} \,
\operatorname{det}\!
{\left[v_{i}^{n-j}-
(v_{i} -a)^{n-j} \: \frac{v_i+a}{v_i}
\prod_{k=1}^n \frac{v_i-u_k}{v_i-u_k-a}\right]}_{1\le i,j \le n}.
\end{multline*}

Our theorems in Section~\ref{sectmain} give rise,
after letting $q=t^\alpha$ and taking the limit $t \rightarrow 1$,
to the following results.
\begin{theo}
Let  $\la=(\la_1,...,\la_{n+1})$ be an arbitrary partition with length
$n+1$. For all $1\le k \le n$ define
$u_k:=\la_k-\la_{n+1}+(n-k)/\alpha$. We have
\begin{equation*}
Q_{(\la_1,\ldots,\la_{n+1})}= \sum_{\ka\in\mathsf{N}^n}
c^{(1/\alpha)}_{\ka_1,\ldots,\ka_n} (u_1,\ldots,u_n;
u_1+{\ka_1},\ldots,u_n+{\ka_n}) \:
Q_{(\la_{n+1}-|\ka|)} \:
Q_{(\la_1+\ka_1,\ldots,\la_n+\ka_n)}.
\end{equation*}
\end{theo}

\begin{theo}
Let  $\la=(1^{m_1},2^{m_2},\ldots,(n+1)^{m_{n+1}})$ be an arbitrary
partition consisting of parts
at most equal to $n+1$. For all $1\le k \le n$ define
$u_k:=\sum_{j=k}^n m_j+(n-k)\alpha$. We have
\begin{multline*}
P_{(1^{m_1},2^{m_2},\ldots,(n+1)^{m_{n+1}})}=
\sum_{\ka\in\mathsf{N}^n}
c^{(\alpha)}_{\ka_1,\ldots,\ka_n} (u_1,\ldots,u_n;
u_1+\ka_1,\ldots,u_n+\ka_n) \\
\times e_{m_{n+1}-|\ka|} \,
P_{(1^{m_1+\ka_1-\ka_2},2^{m_2+\ka_2-\ka_3},\ldots,
{(n-1)}^{m_{n-1}+\ka_{n-1}-\ka_n},n^{m_n+m_{n+1}+\ka_n})}.
\end{multline*}
\end{theo}

As in Section~\ref{sectanal} these formulas generate explicit analytic 
developments
for Jack polynomials. 

\medskip
{\footnotesize \noindent
\textbf{Acknowledgements :}
The first author is pleased to thank Alain Lascoux for friendly advices.}

\bigskip
\noindent{\begin{tabular}{p{9cm}p{11.5cm}}
    Michel Lassalle&Michael Schlosser\\ 
    {\footnotesize Centre National de la Recherche Scientifique}
        &{\footnotesize Institut f\"{u}r Mathematik}\\
    {\footnotesize Institut Gaspard Monge}
        &{\footnotesize Universit\"{a}t Wien}\\
    {\footnotesize 77454 Marne-la-Vall\'ee Cedex, France}
        &{\footnotesize A-1090 Wien, Autriche}\\
    {}&{}\\ 
    {\footnotesize {e-mail: lassalle @ univ-mlv.fr}}&
    {\footnotesize {e-mail: schlosse @ ap.univie.ac.at}}\\
    {\footnotesize {http://www-igm.univ-mlv.fr/
    {\textasciitilde}lassalle/index.html}}
        &{\footnotesize 
{http://www.mat.univie.ac.at/{\textasciitilde}schlosse/}}\\ 
    \end{tabular}}
    
\end{document}